\documentclass[twoside,11pt]{article}
\usepackage{jmlr2e} 
\usepackage{amsmath, amsfonts, amssymb, epsfig, pstricks,pst-node,pst-tree,graphicx,./aoslay}

\usepackage{flushend}
\newcommand{\Ab}[1]{A_{{}_{#1}}}
\newcommand{\Exp}{\mathbb E}

\newcommand{\dist}{{\cal J}}
\newcommand{\eps}{\epsilon_{{}_{p'}}}
\newcommand{\R}{R}

\renewcommand{\footnoterule}{\kern -3pt \hrule \kern 2.6pt}

\ShortHeadings{Insurability}{Santhanam and Anantharam}

\firstpageno{1} 

\begin{document}

\title{Agnostic insurability of model classes}
\author{\name Narayana Santhanam \email nsanthan@hawaii.edu\\
\addr Dept of Electrical Engineering\\
University of Hawaii at Manoa\\
Honolulu, HI 96822
\AND
\name Venkat Anantharam \email ananth@eecs.berkeley.edu\\
\addr Dept of EECS\\
University of California, Berkeley\\
Berkeley, CA 94720}
\editor{}
\maketitle

\begin{abstract}
Motivated by problems in insurance, our task is to predict finite
upper bounds on a future draw from an unknown distribution $p$ over
the set of natural numbers. We can only use past observations generated
independently and identically distributed according to $p$. While $p$
is unknown, it is known to belong to a given collection $\cP$ of
probability distributions on the natural numbers. 

The support of the distributions $p \in \cP$ may be unbounded, and the prediction game goes on for \emph{infinitely}
many draws. We are allowed to make observations without predicting upper bounds for some time. But we must, with
probability $1$, start and then continue to predict upper bounds after a finite time irrespective of which $p \in \cP$
governs the data.

If it is possible, without knowledge of $p$ and for any prescribed
confidence however close to $1$, to come up with a sequence of upper
bounds that is never violated over an infinite time window with
confidence at least as big as prescribed, we say the model
class $\cP$ is \emph{insurable}.  

We completely characterize the insurability of any class $\cP$ of
distributions over natural numbers by means of a condition on how the
neighborhoods of distributions in $\cP$ should be, one that is both
necessary and sufficient.
\end{abstract}


\section{Introduction}	\label{s.intro}

Insurance is a means of managing risk by transfering a potential
sequence of losses to an \emph{insurer} for a price paid on a 
regular basis, the \emph{premium}. The insurer
attempts to break even by balancing the possible loss that may be
suffered by a few with the guaranteed premiums of many. We aim to
study the fundamentals of this problem when the losses can be
unbounded and a precise model for the probability distribution of the
aggregate loss in each period either does not exist or is infeasible
to get.

A systematic, theoretical, as opposed to empirical, study of insurance
goes back to 1903 when Filip Lundberg (see~\cite{EM01}) defined a natural
probabilistic setting as part of his thesis. In particular, Lundberg
formulated a collective risk problem pooling together the risk of all
the insured parties into a single entity, which we call the insured.
Typically, studies of insurance derived from the approach
in~\cite{EM01} depend on working with specific models for the loss
distribution, e.g.  compound Poisson models, after which questions of
interest in practice, such as the relation between the
size of the premiums charged and the probability of the insurer going
bankrupt, can be analyzed.  A rather comprehensive theory of insurance
along these lines has evolved in~\cite{Cra69} and more recently
in~\cite{Asm}. They incorporate several model classes for the
distribution of the losses over time other than compound Poisson
processes, including some heavy tailed distribution classes.

\paragraph{}
We depart from the existing literature on insurance in two important respects.

\paragraph{No upper bound on loss} The first departure relates to the
practice among insurers to limit payments to a predetermined ceiling,
even if the loss suffered by the insured exceeds this ceiling. In both
the insurance industry and the legal regulatory framework surrounding
it, this is assumed to be common sense. But is it always necessary to
impose such ceilings? Moreover, in scenarios such as reinsurance, a
ceiling on compensation is not only undesirable, but may also limit
the very utility of the business. As we will see, we may be able to
handle scenarios where the loss can be unbounded.

\paragraph{Universal approach} The second aspect of our
approach arises from our motivation to deal with several new settings
for which some sort of insurance is desirable, but where insurers are
hesitant to enter the market due to lack of sufficient data.  Examples
of such settings include insuring against network outages or attacks
against future smart grids, where the cascade effect of outages or
attacks could be catastrophic.  In these settings, it is not clear
today what should constitute a reasonable risk model because of the
absence of usable information about what might cause the outages or
motivate the attacks.

We address the second issue by working with a \emph{class} of models, 
\ie a set of probability laws over loss sequences that adheres to any 
assumptions the insurer may want to make or any information it may 
already have. In this paper we will only consider loss models that are 
independent and identically distributed (\iid) from period to period, 
so we can equivalently think of a model class as defined in terms of 
its one dimensional marginals. 

As an example, we may want to consider the set of all finite moment
probability distributions over the nonnegative integers as our class
of possible models for the loss distribution in each period. Now, we
ask the question: what classes of models are the ones on which the
insurer can learn from observations and set premiums so as to remain
solvent?  In this paper, we completely answer this question by giving
a necessary and sufficient condition that characterizes what classes
of models lend themselves to this insurance task.

This is very reminiscent of the universal compression/estimation/prediction 
approaches (see~\cite{Sht87,Fit72,Ris84,Rya08})---we will have more to
say on this shortly. There is also extensive work regarding learning 
from experts that has a related flavor, see~\cite{CL} for a survey.

\paragraph{Formulation} Formally, we adopt the collective risk
approach, namely, we abstract the problem to include just two agents,
the insurer and the insured.  Losses incurred by the insured are
considered to form a discrete time sequence of random variables, with
the sequence of losses denoted by $\sets{X_i,~i\ge1}$, and we assume
that $X_i\in\naturals$ for all $i\ge 1$, where $\naturals$ denotes the
set of natural numbers, $\sets{0,1,2,\ldots}$.

A \emph{model class} $\cP^{\infty}$ is a collection of measures on
infinite length loss sequences, and is to be thought of as the set of
all potential probability laws governing the loss sequence. Each
element of $\cP^{\infty}$ is a \emph{model} for the sequence of
losses. Any prior knowledge on the structure of the problem is
accounted for in the definition of $\cP^\infty$. We focus on measures
corresponding to \iid samples, i.e. each member of $\cP^{\infty}$
induces marginals that are product distributions. We denote by $\cP$ the set of
distributions on $\naturals$ obtained as one dimensional marginals of
$\cP^{\infty}$. Since there is no risk of confusion, we will also
refer to the distributions in $\cP$ as models and to $\cP$ as the
model class.

The actual model in $\cP$ governing the law of the loss in each period
remains unknown to the insurer.  We assume no ceiling on the loss, and
require the insurer to compensate the insured in full for the loss in
each period at the end of that period.  The insurer is assumed to
start with some initial capital $\Pi_0 \in \reals^+$, a nonnegative
real number.  The insurer then sets a sequence of premiums based on
the past losses---at time $i$, the insurer collects a premium
$\Pi(X_1^{i-1})$ at the beginning of the period, and pays out full
compensation for loss $X_i$ at the end of the period. If the built up
capital till step $i$ (including $\Pi(X_1^{i-1})$, and after having
paid out all past losses) is less than $X_i$, the insurer is said to
be \emph{bankrupted}.

Given a class $\cP^{\infty}$ of loss models, we ask if for every
prescribed upper bound $\eta > 0$ on the probability of bankruptcy,
the insurer can set (finite) premiums at every time step based only on
the loss sequence observed thus far and with no further knowledge of
which law $p \in \cP^{\infty}$ governs the loss sequence, while
simultaneously ensuring that the insurer remains solvent with
probability bigger than $1 - \eta$ under $p$ irrespective of which $p
\in \cP^{\infty}$ is in effect. If the probability of the insurer ever
going bankrupt over an infinite time window can be made arbitrarily
small in this sense, the class of \iid loss measures $\cP^\infty$ is
said to be \emph{insurable}.

A couple of clarifications are in order here.  First, to make the
problem non-trivial, we allow the insurer to observe the loss sequence
for some arbitrary finite length of time without having to provide
compensations.  We require that the insurer has to eventually provide
insurance with probability 1 no matter which $p \in \cP^{\infty}$ is
in effect. The insurer cannot quit providing insurance once it has
entered into the insurance contract with the insured.  Premiums set
before the entry time can be thought of as being $0$ and the question
of bankruptcy only arises after the insurer has entered into the
contract.  Secondly, at this point of research, we do not concern
ourselves with incentive compatibility issues on the part of the
insured and assume that the insured will accept the contract once the
insurer has entered, agreeing to pay the premiums as set by the
insurer.

It turns out that the fact that the capital available to the insurer
at any time is built up from past premiums does not play any role in
whether a model class is insurable or not. In fact, the problem is
basically one of finding a sequence of finite upper bounds
$\Phi(X_1^{i-1})$ on the loss $X_i$ for all $i \ge 1$. We refer to the
sequence $\{ \Phi(X_1^{i-1}),~i \ge 1\}$ as the \emph{loss dominating
  sequence} and call $\Phi(X_1^{i-1})$ the \emph{loss-dominant} at
step $i$.

The notion of insurability of a model class $\cP$ comes down to
whether for each $\eta > 0$ there is a way of choosing the loss
dominants such that the probability of the loss $X_i$ ever
exceeding the loss dominant $\Phi(X_1^{i-1})$ is smaller than $\eta$
irrespective of which model $p$ in the model class $\cP^{\infty}$ is
in effect. Here again we allow some initial finite number of periods
for which the loss dominant can be set to $\infty$, but it must become
finite with probability $1$ under each $p \in \cP^{\infty}$ and stay
finite from that point onwards.

\paragraph{}
It will be interesting to examine this formulation in the broader
context of pointwise convergent algorithms, in particular universal
compression algorithms.

\paragraph{Pointwise convergence}
Theoretically, the flexibility we have permitted regarding when to
start proposing finite loss dominants allows us to categorize the
insurance problem formulated above as one that admits what we call
\emph{useful} pointwise convergent estimators, even when uniformly
convergent estimates are impossible. Roughly speaking, the insurance
problem can be thought of requiring estimation of all the percentiles
of an unknown distribution from $\cP$, using only \iid draws generated
from the distribution. However, as the sample size increases, the estimate of any
given percentile need not converge to the true value (according to
some predefined metric) uniformly over the entire class $\cP$.

In general, estimators whose rate of convergence cannot be bounded by
parameters that are known \emph{a-priori} or observed from the sample
are often frowned upon by practitioners. This is because even if we
know that such an estimator is consistent, for a given sample there
may be no way of telling of how good or bad the estimate is. 

This poses a connundrum, since when dealing with large alphabets or
high dimensions, it is sometimes too restrictive to only deal with
model classes or problem formulations that admit uniformly convergent
estimators---those that converge to the true values at a rate that can
be bounded uniformly over the model class as the sample size increases
to infinity.

What if we are forced to work with a model class which is sufficiently
complex that uniformly convergent estimators are impossible? We can
still salvage the situation if for any given finite sample, we had
some way to tell if the estimate was doing well or not relative to the
true unknown model, perhaps by looking at the sample at hand.

This is what our insurance formulation capitalizes on as
well. Requiring only classes $\cP$ of distributions that have
uniformly convergent estimators for percentiles of a distribution is
too restrictive when the support of $\cP$ is not bounded. We want to
deal with broader classes of models, and the flexibility about the
start point for prediction allows us to consider significantly richer
model classes. For other kinds of such ``useful'' pointwise estimation,
particularly in relation to Markov processes, see~\cite{ARS13:arxiv}.

\paragraph{Universal compression} The approach we take is not
unconnected with universal compression literature, as well as learning 
formulations involving regret. The point of departure is that our approaches 
are interesting precisely in cases where the strong notions of (worst-case 
or average-case) redundancy fail. Namely, classes of distributions whose
redundancy is not finite.

The closest universal compression formulation to our notion of insurability 
here is in the idea of weak universal compression in~\cite{Kie78}. However, 
weak compression does not include the aspect of determining from the data at
hand when a compressor is doing well---a crucial part of our problem.

Even so, the relation between this problem and weak compression is far 
from obvious. Perhaps surprisingly, it is completely possible that we can 
find (see~\cite{SA12:spcom}) classes of models that can be compressed weakly 
but are not insurable, as well as classes of models that are insurable, 
but cannot be compressed weakly. 

\paragraph{}
However, there is one insight that we conjecture can be generalized
beyond insurability to all problems with the flavor of useful
pointwise convergence of estimates---local complexity as opposed to
global complexity of model classes. Insurability of model classes does
not depend on global complexity measures of model classes---such the
redundancy of model classes or the Rademacher complexity. Instead,
insurability is related to how the local neighborhoods look like, as
we will see in Section~\ref{s:rsl}.

\paragraph{Results} For a model class to be insurable, roughly speaking, close
distributions must have comparable percentiles. Distributions in the
model class that, in every neighborhood, have some other distribution
with arbitrarily different percentiles are said to be
\emph{deceptive}. In Section~\ref{s:rsl}, we define what it means for
distributions to be close, and what it means for distributions to have
comparable percentiles. In Section~\ref{s:ex}, we provide several
examples of insurable and non-insurable model classes.  Our main
result is Theorem~\ref{thm:ncssff} of Section~\ref{s:rsl}, which states that that
$\cP^\infty$ is insurable iff it has no deceptive distributions. We
prove this theorem in Sections~\ref{s:ncs} and~\ref{s:sff}.

\ignore{Let $\naturals^{*}$ be the collection of all finite length
  strings of natural numbers. The insurer's \emph{scheme} $\Phi$ is a
  mapping from $\naturals^*\to\reals^+$, and is interpreted as the
  premium demanded by the insurer from the insured after a loss
  sequence is observed. The insurer can observe the loss for a time
  prior to entering the insurance game.  However, we require the
  insurer enters the game with probability 1 no matter what loss model
  is in force, and the insurer cannot quit once entered. In
  Theorems~\ref{thm:ncs} and~\ref{thm:sff}, we determine a necessary
  and sufficient condition on $\cP$ for insurability.}

\section{Precise formulation of the problem}
\label{s:prob}
We model the loss at each time by a random variable taking values in
$\naturals=\sets{0,1,\ldots}$.  Denote the sequence of losses by $X_1,
X_2\ldots$ where $X_i\in \naturals$.  Let $\naturals^*$ be the set of
all finite length sequences from $\naturals$, including the empty
sequence. We will write $x^n$ for the sequence $x_1\upto x_n$.  Where
it appears, $x^0$ denotes the empty sequence. A loss distribution is a
probability distribution on $\naturals$. Let $\cP$ be a set of loss
distributions. $\cP^{\infty}$ is the collection of \iid measures over
infinite sequences of symbols from $\naturals$ such that the set of
one dimensional marginals over $\naturals$ they induce is $\cP$.

We write $\reals^+$ for the set of nonnegative real numbers and use
$:=$ for equality by definition.

Consider an insurer with an \emph{initial capital} $\Pi_0 \in
\reals^+$.  An \emph{insurance scheme} for $\cP$ is comprised of a
pair $(\tau, \Pi)$.

Here $\tau ~:~ \naturals^* \mapsto \{0, 1\}$ satisfies $\tau(x_1,
\ldots, x_n) = 1 \Longrightarrow \tau(x_1, \ldots, x_{n+1}) = 1$ for
all $x^n$ and also $p( \sup_n \tau(X^n) = 1) = 1$ for all $p \in
\cP^{\infty}$.  $\tau$ should be thought of as defining an \emph{entry
  time} for the insurer with the property that once the insurer has
entered it stays entered and that the insurer enters with probability
$1$ irrespective of which $p \in \cP^{\infty}$ is in effect. Here we
say the insurer enters after seeing the sequence $x^n \in \naturals^*$
(possibly the empty sequence) if $\tau(x^n) =1$. The other ingredient
of an insurance scheme is the \emph{premium setting scheme}
$\Pi:\naturals^*\to\reals^+$, satisfying $\Pi(x^n) = 0$ if $\tau(x^n)
=0$, with $\Pi(x^n)$ being interpreted as the premium demanded by the
insurer from the insured after the loss sequence $x^n \in \naturals^*$
is observed.

Let $1(\cdot)$ denote the indicator function of its argument.  The
event that the insurer goes bankrupt is the event that
\[
\Pi_0 + \sum_{i=1}^n (\Pi(X^{i-1}) - X_i)1(\tau(X^{i-1}) =1) < 0 \mbox{ for some $n \ge 1$}~.
\]
In words, this is the event that in some period $n \ge 1$ after the
insurer has entered, the loss $X_n$ incurred by the insured exceeds
the built up capital of the insurer, namely the sum of its initial
capital and all the premiums it has collected after it has entered
(including the currenly charged premium $\Pi(X^{n-1})$) less all the
losses paid out so far.


\bDefinition A class $\cP^{\infty}$ of laws on loss sequences is
called insurable by an insurer with initial capital $\Pi_0 \in
\reals^+$ if $\forall$ $\eta>0$, there exists an insurance scheme
$(\tau, \Pi)$ such that $\forall$ $p\in\cP^{\infty}$,
\[
p\Paren{(\tau, \Pi) \text{ goes bankrupt }}<\eta~.
\] 
\eDefinitionp

We should remark that despite the apparent role of the initial capital
of the insurer in this definition, it plays no role from a
mathematical point of view. To see this note first that if a model
class $\cP^{\infty}$ is insurable by an insurer with capital $\Pi_0$
it is clearly insurable by all insurers with initial captial at least
$\Pi_0$, since such an insurer can use the same entry time and premium
setting scheme as the insurer with initial capital $\Pi_0$.  On the
other hand, an insurer with initial capital less than $\Pi_0$ can use
the same entry time as an insurer with initial capital $\Pi_0$ and
simply charge an additional premium at the time of entry which in
effect builds up its initial capital to $\Pi_0$, and then proceed with
the same premium setting scheme as that used by the insurer with
initial captial $\Pi_0$. This feature is an artifact of the complete
flexibility we give the insurer in setting premiums; for more on this
see the concluding remarks in Section \ref{s.concrem}.

As indicated in the introductory Section \ref{s.intro}, 
we will first show that whether a model class
of loss distributions is insurable is equivalent to whether we can find suitable 
loss domination sequences for the sequence of losses. 
We next make this connection and the associated terminology precise.

\bDefinition A \emph{loss domination scheme} for $\cP$ is a mapping
$\Phi ~:~ \naturals^* \mapsto \reals^+ \cup \{ \infty \}$, where for 
$x^n\in \naturals^*$, we interpret
$\Phi(x^n)$ as an estimated upper bound on $x_{n+1}$. 
We call
$\{ \Phi(X^{i-1}),~i \ge 1\}$ the loss-domination sequence and $\Phi(X^{i-1})$ 
the loss-dominant at step $i$. 
We require for all $x^n \in \naturals^*$ that
\[
\Phi(x_1, \ldots, x_n) < \infty \Longrightarrow 
\Phi(x_1, \ldots, x_{n+1}) < \infty
\] 
and also that for all $p \in \cP^{\infty}$, 
\[
p( \inf_{n \ge 1} \Phi(X^n) < \infty) = 1.\eqed
\] 
\eDefinitionp
We think of $\Phi(x^n) = \infty$ as saying that the scheme has not yet committed to 
proposing finite loss dominants after having seen the sequence $x^n$, while if 
$\Phi(x^n) < \infty$ it has. Once the scheme commits to proposing finite loss dominants
it has to continue to propose finite loss dominants from that point onwards.
Further, with probability $1$ under every $p \in \cP^{\infty}$, the scheme has to eventually 
start proposing finite loss dominants. 

\bDefinition
Given our motivation from the insurance problem,
we will say the loss domination scheme $\Phi$ goes \emph{bankrupt}
if $\Phi(X^{n-1}) < X_n$ for some $n \ge 1$.
\eDefinition

The connection between the insurance problem and the problem of selecting loss
dominants can now be made precise as follows.

\bObservation 
Let $\cP^{\infty}$ be a model class and $\eta > 0$. Let $\Pi_0 \in \reals^+$. 
An insurer with initial capital $\Pi_0$ can find an insurance scheme $(\tau, \Pi)$
such that the probability of remaining solvent is bigger than $1 - \eta$ irrespective of which
$p \in \cP^{\infty}$ is in effect if and only if there is a loss domination scheme 
$\Phi$ such that the probability of it going bankrupt is less than $\eta$ irrespective
of which $p \in \cP^{\infty}$ is in effect.

\ignore{
\bObservation Let $P$ be an arbitrary measure on $\naturals^*$ and
let $\eta>0$.
There exists premiums $\Pi$ such that
\begin{equation}
\label{eq:pi}
P\Paren{ \exists i:
\sum_{j=1}^{i}
(\Pi(X_1^{j-1})-X_j){\bf 1}(\text{ $\Pi$ has entered before step $j$ }) 
<0 } < \eta,
\end{equation}
iff there is a scheme $\Phi$ that, for every loss sequence enters the game
at the same points $\Pi$ does, and satisfies 
\begin{equation}
\label{eq:phi}
P\Paren{ \exists i:
(\Pi(X^{i-1})-X_i){\bf 1}(\text{ $\Pi$ has entered before step $i$ }) 
<0 } < \eta.
\end{equation}
satisfying the bankruptcy condition, note that
in any round $i$, the scheme
\[
\Phi(X^{i-1})=\sum_{j=1}^{i}\Pi(X_1^{j-1}){\bf 1}(\text{ $\Pi$ has entered before step $j$ }) 
\] 
exceeds the sum of all built-up profits and the premium
$\Pi(X^{i-1})$ charged at step $i$. Given any sequence of losses, if
$\Pi$ survives bankruptcy till step $i$, so does $\Phi$.  Therefore,
if $\Pi$ satisfies~\eqref{eq:pi}, $\Phi$ satisfies~\eqref{eq:phi}.
Now given any scheme $\Phi'$ satisfying
\[
P\Paren{ \exists i:
(\Phi'(X^{i-1})-X_i){\bf 1}(\text{ $\Phi'$ has entered before step $i$ }) 
<0 } < \eta.
\] 
let us set premiums as $\Pi'=\Phi'$. Again given any sequence of losses, 
if $\Phi'$ survives till step $i$, $\Pi'$ does so without even requiring
the built up profits. Hence the above observation.~\hfill$\Box$
}

\Proof
Given an insurance scheme $(\tau,\Pi)$ consider the loss domination scheme $\Phi$ that
has $\Phi(x^n) := \infty$ iff $\tau(x^n) = 0$ and 
\[
\Phi(X^{n-1}) := \Pi_0 + \sum_{i=1}^{n-1} (\Pi(X^{i-1}) - X_i)1(\tau(X^{i-1}) =1)\, +\, \Pi(X^{n-1})~,
\]
if $\tau(X^n) = 1$. Since $\tau$ enters (become equal to $1$) with probability $1$ 
under each $p \in \cP^{\infty}$ and stays equal to $1$ once it has become $1$,
$\Phi$ becomes finite with probability $1$ under each $p \in \cP^{\infty}$ and stays
finite once it has become finite. Thus $\Phi$ is indeed a loss domination scheme.
It is straightforward to check that if the insurance scheme $(\tau,\Pi)$ 
stays solvent with probability bigger than $1 - \eta$ irrespective of which
$p \in \cP^{\infty}$ is in effect then the loss domination scheme $\Phi$ becomes bankrupt
with probability less than $\eta$ irrespective of which $p \in \cP^{\infty}$ is in effect.

Conversely, given a loss domination scheme $\Phi$ define the insurance scheme 
$(\tau,\Pi)$ by setting $\tau(x^n) := 0$ iff $\Phi(X^n)  = \infty$ (and 
$\tau(x^n) := 1$ iff $\Phi(x^n) < \infty$) and defining $\Pi(x^n) := 0$ if
$\Phi(x^n) = \infty$ and $\Pi(x^n) := \Phi(x^n)$ if $\Phi(x^n) < \infty$. 

One sees that $\tau$ as defined becomes $1$ with probability $1$ under each 
$p \in \cP^{\infty}$ and stays equal to $1$ once it becomes $1$. Further, the premiums
set at each time are finite and equal to $0$ till the entry time. Thus $(\tau, \Pi)$ 
as defined is indeed an insurance scheme.

It is straightforward to check if $\Phi$ becomes bankrupt with probability less than $\eta$ irrespective of which $p \in
\cP^{\infty}$ is in effect, then $(\tau,\Pi)$ stays solvent with probability bigger than $1 - \eta$ irrespective of
which $p \in \cP^{\infty}$ is in effect.  Hence the above observation.~\hfill$\Box$ \eObservation

We may therefore conclude that a model class $\cP^{\infty}$ is insurable iff 
for all $\eta > 0$ there is a loss domination scheme $\Phi$ such that the probability of
going bankrupt under $\Phi$ is less than $\eta$ irrespective of which $p \in \cP^{\infty}$
is in effect. In the rest of the paper we will therefore 
focus mainly on whether the model class $\cP^{\infty}$ is such that for every 
$\eta > 0$ a loss domination
sequence $\Phi$ exists with its probability of bankruptcy being less than $\eta$
irrespective of which model in the model class governs the sequence of losses.

In Theorem~\ref{thm:ncssff}, we provide a condition
on $\cP$ that is both necessary and sufficient for insurability.

\section{Statement of the main result}
\label{s:rsl}
We go through a few technical points before spelling out the 
results in detail in~\ref{ss:rsl}.
\subsection{Close distributions}
Insurability of $\cP^\infty$ depends on the neighborhoods of the
probability distributions among its one dimensional marginals $\cP$. The
relevant ``distance'' between distributions in $\cP$ that decides the
neighborhoods is
\[
\dist(p,q)
:=
D\Paren{p||\frac{p+q}2}
+
D\Paren{q||\frac{p+q}2}.
\]
Here $D(p||q)$ denotes the relative entropy of $p$ with respect to $q$, where
$p$ and $q$ are probability distributions on $\naturals$, defined by 
\[
D(p||q) := \sum_{y \in \naturals} p(y) \log \frac{p(y)}{q(y)}~.
\]
The logarithm is assumed to be taken to base $2$ (we use $\ln$ for the logarithm
to the natural base).

\subsection{Cumulative distribution function}
\label{s:cdf}

Since we would like to discuss percentiles, it is convenient to use a non-standard
definition for the cumulative distribution function of a probability distribution on 
$\naturals$.

For our purposes, the cumulative distribution function of any probability distribution
$p$ on $\naturals$ is a function from 
$\reals^+ \cup \{\infty\} \to [0,1]$, and will be denoted by $F_p$.
We obtain $F_p$ by first defining $F_p$ on points in the support of
$p$.
We define $F_p$ for all other 
nonnegative real numbers
by linearly interpolating between the values in the support of
$p$. Finally, $F_p(\infty) := 1$.

Let $F_p^{-1} ~:~ [0,1] \mapsto \reals^+ \cup \{\infty\}$ denote the 
inverse function of $F_p$. Then $F_p^{-1}(x) =0$ for all $0\le x< F_p(0)$.
If $p$ has infinite support then $F_p^{-1}(1) = \infty$, else
$F^{-1}_p(1)$ is the smallest natural number $y$ such that
$F_p(y)=1$.

Two 
simple and useful observations can now be made.
Consider a probability distribution $p$ with support $\cA\subset
\naturals$.  For $\delta>0$, let ($T$ for tail)
\[
T_{p,\delta} := \sets{y\in \cA: y\ge F^{-1}(1-\delta)},
\] 
and let ($H$ for head)
\[
H_{p,\delta} := \sets{y\in \cA: y\le 2F^{-1}(1-\delta/2)}.
\] 
It is easy to see that 
\begin{equation}
\label{eq:tail}
p(T_{p,\delta})> \delta
\end{equation} and that
\begin{equation}
\label{eq:head}
p(H_{p,\delta})>1-\delta.
\end{equation}
Suppose that for some $\delta >0$ we have $F_p^{-1}(1-\delta)>0$
and the loss-dominant at the beginning of period $i \ge 1$ happens to be
set to $F_p^{-1}(1-\delta)$, then the probability under $p$ of the loss in period $i$
exceeding the loss-dominant is bigger than $\delta$.  If the loss-dominant 
at the beginning of period $i$ happens to be set to
$2F_p^{-1}(1-\delta/2)$, then the probability that the loss in period $i$ exceeds the
loss-dominant is less than $\delta$. We will use these observations in the proofs
to follow.

\subsection{Necessary and sufficient conditions for insurability}
\label{ss:rsl}
Existence of close distributions with very different quantiles is what
kills insurability. A loss domination scheme could be ``deceived'' by some process
$p\in\cP^\infty$ into setting low loss-dominants, while a close enough
distribution hits the scheme with too high a loss. The conditions for insurability
of $\cP^\infty$ are phrased in terms of the set of its one dimensional marginals,
$\cP$.

Formally, a probability distribution $p$ in $\cP$ is called \emph{deceptive} if
$\forall$ $\epsilon>0$, $\exists$
$\delta>0$ such that that no matter what $f(\delta)\in\reals^+$ is chosen,
$\exists$ a (bad) distribution $q\in\cP$ such that 
\[
\dist(p,q) < \epsilon
\] 
and 
\[
F_q^{-1}(1-\delta)> f(\delta).
\]
In the above definition, $f(\delta)$ is simply an arbitrary
nonnegative real number. However, it is useful to think of this number as the
evaluation of a function $f:(0,1)\to\reals$ at $\delta$.
Equivalently, a distribution $p$ in $\cP$ is not \emph{deceptive} if
$\exists$ $\epsilon_p>0$, such that $\forall$ $\delta>0$,
$\exists$ $f(\delta)\in\reals$, such that all distributions $q\in\cP$
with
\[
\dist(p,q) < \epsilon
\] 
satisfy
\[
F_q^{-1}(1-\delta) \le f(\delta).
\]
\ignore{If $p$ is not deceptive, then for all $\delta_p>0$, $\exists$ $f_p:
\reals\to \reals$, $\epsilon_{p}>0$, such that all $q$ satisfying
$\dist(p,q)<\epsilon_p$ are $f_p$-matched to $p$ on tail
$\delta_p$.}

Our main theorem is the following, which we prove in 
Sections~\ref{s:ncs} and~\ref{s:sff}.
 
\bTheorem
\label{thm:ncssff} 
$\cP^{\infty}$ is insurable, iff no $p\in\cP$ is deceptive. 
\eTheorem


\section{Examples}
\label{s:ex}

Consider $\cU$, the collection of all uniform distributions over a
finite contiguous support of the form $\sets{m\upto M}$, with $m \le M$ being
arbitrary nonnegative integers. Let the losses come as \iid samples from one of the
distributions in $\cU$---call the resulting model class $\cU^{\infty}$.
\bExample $\cU^\infty$ is insurable.  
\Proof If the threshold
probability of ruin is $\eta$, choose the loss-domination scheme $\Phi$ as follows. For
all sequences $x^n$ with $n \le \log \frac1{\eta}+1$ set
$\Phi(x^n)=\infty$.  For all sequences $x^n$ with $n > \log
\frac1{\eta}+1$, the loss-dominant $\Phi(x^n)$ is set to be 
twice the largest loss observed thus
far. It is easy to see that this scheme is bankrupted with probability
less than $\eta$ irrespective of which $p \in \cU^{\infty}$ is in effect.  \eExample

Consider the set $\cN^{\infty}$ of all \iid processes such that
the one dimensional marginals have finite moment. Namely, $\forall
p\in\cN^{\infty}$, $\Exp_pX_1<\infty$.

\bExample\label{ex:bnd} $\cN^{\infty}$ is not insurable.  

\Proof Note
that the loss process that puts probability 1 on the all zero
sequence exists in $\cN^{\infty}$, since it corresponds to
the 
one dimensional marginal loss distribution 
that produces loss $0$ in each period.
Since every loss domination scheme
enters with probability 1 no matter which $p\in\cN^{\infty}$ is in
force, every loss domination scheme must enter after seeing some finite number
of zeros.  Fix any loss domination scheme $\Phi$. 
Suppose the scheme starts to set finite dominants after seeing $N$
losses of size 0. To show that $\cN^{\infty}$ is not insurable, we
show that $\exists \eta>0$ 
and $\exists
p\in\cN^{\infty}$ such that
\[
p(\text{ $\Phi$ goes bankrupt })\ge \eta.
\] 
Fix $\delta=1-\eta$. 
Let $\epsilon$ be small enough that 
\[
(1-\epsilon)^N> 1-\delta/2, 
\]
and let $M$ be a number large enough that 
\[
(1-\epsilon)^M<\delta/2.
\]
Note that since $1-\delta/2\ge \delta/2$, we have $N< M$. 
Let $L$ be greater than any of loss-dominants set by $\Phi$
for the sequences $0^N, 0^{N+1},\ldots 0^M$. 
Let $p\in\cN^{\infty}$ satisfy, for all $i$,
\[
p(X_i)=
\begin{cases}
1-\epsilon &\text{ if } X_i=0\\
\epsilon & \text{ if } X_i=L.
\end{cases}
\]
For the \iid loss process having the law $p$, the insurer is bankrupted on all sequences that
contain loss $L$ in between the $N$-th and $M$-th steps. These sequences,
$0^NL, 0^{N+1}L\upto 0^{M-1}L$, have respective probabilities (under $p$) 
\[(1-\epsilon)^N\epsilon, (1-\epsilon)^{N+1}\epsilon,\ldots, (1-\epsilon)^{M-1},
\] 
and they also form a prefix free set. Therefore, summing up the
geometric series and using the assumptions on $\epsilon$ above,
\[
p(\text{ $\Phi$ is bankrupted })\ge
(1-\epsilon)^N-(1-\epsilon)^M
\ge
1-\delta/2 -\delta/2=\eta.\eqed
\]  
\eExamplep 
One can actually directly verify that every distribution in
$\cN^{\infty}$ is deceptive.

Consider the collection of all \iid loss distributions with monotone one dimensional marginals. A monotone probability
distribution $p$ on $\naturals$ is one that satisfies $p(y+1) \le p(y)$ for all $y \in \naturals$.  Let $\cM^{\infty}$
be the set of all \iid loss processes, with one dimensional marginal distribution from $\cM$, the collection of all
monotone probability distributions over $\naturals$.

Again, it is easily shown that every distribution in $\cM$ is deceptive.  It follows from Theorem~\ref{thm:ncssff} that
\bExample
\label{ex:mnt}
$\cM^{\infty}$ is not insurable.
\eExample
 
Now for $h>0$, we consider the set $\cM_h\subset\cM$ of all
monotone distributions over $\naturals$ whose entropy is upper bounded
by $h$. Let $\cM_h^\infty$ be the set of all \iid loss processes with
one dimensional marginals from $\cM_h$. Then

\bExample
\label{ex:mnth}
$\cM_h^\infty$ is insurable.
\Proof From Markov inequality, if $p\in\cM_h$ and $X\sim p$, 
\[
p( X> M ) = p(\log X > \log M) < \frac{E_p\log X}{\log M} \le \frac{E_p\log \frac1{p(X)}}{\log M} \le
\frac{h}{\log M}.
\]
To see the second inequality above, note that $p$ is monotone therefore for any number $i$, $p(i)\le \frac1i$. Therefore, for
all $p\in\cM_h$, 
\[
F_p^{-1}(1-\delta) \le 2^{\frac{H}{\delta}}.
\]
Thus no $p\in\cM_h$ is deceptive, and $\cM_h^\infty$ is insurable.
\eExample

\section{Necessary condition for insurability}
\label{s:ncs}

\ignore{
Recall that with the unconventional definition of
cumulative distribution functions in Section~\ref{s:cdf}, if for a
sequence $x^n$ we have $F_q^{-1}(1-\delta) > \Phi(x^n)$, then 
the loss domination scheme $\Phi$ will be
bankrupted with probability at least $\delta$ in the next step.
Recall that the model class $\cP^{\infty}$ is a set of \iid measures over infinite sequences 
from $\naturals$, and that $\cP$ denotes the collection of one dimensional marginals of
members of $\cP^{\infty}$.
}

In this section we prove one direction of Theorem \ref{thm:ncssff}, as stated next.

\bTheorem
\label{thm:ncs} 
If $\cP^{\infty}$ is insurable, then no $p\in\cP$ is deceptive. 
\Proof 
To keep notation simple, we will denote by $p$ (or $q$) both a measure
in $\cP^{\infty}$ as well as the corresponding one dimensional marginal distribution,
which is a member of $\cP$. The context will clarify which of the two is meant.
We prove the contrapositive of the theorem: if some $p\in\cP$
is deceptive, then $\cP^{\infty}$ is not insurable.  

\ignore{Pick $0\le \alpha< h^{-1}(\half)$ where $h(x)$ is the binary entropy
function defined for $0\le x\le 1$ by
\[
h(x) :=-x\log x -(1-x)\log (1-x),
\]
where the logarithm is to base $2$.}
Pick $\alpha>0$ and fix $0< \eta<(1-\alpha -\frac2N)\Paren{1-\frac1e}$.
Suppose $p\in\cP$ is deceptive.  We prove that $\cP^\infty$ is not
insurable by finding for each loss domination scheme $\Phi$, a
probability distribution $q \in \cP$ close to $p$ such that
\[
q(\text{ $\Phi$ goes bankrupt })\ge \eta.
\] 
The basic idea is that because $\Phi$ has to enter with probability 1
under $p$, it would have been forced to set premiums that are too low
for $q$.
\ignore{To show the contraposition, let $p\in\cP$ denote a deceptive
distribution.  Namely, for all $\epsilon>0$, $\exists\delta>0$ such
that $\forall f(\delta) \in \reals$, $\exists$ $q$
satisfying $\dist(p,q)<\epsilon$ and
\[
F_q^{-1}(1-\delta)> f(\delta).
\]}

\ignore{Note that for the deceptive $p$, $\exists$ $\forall
f_p:\reals\to\reals$ and $\epsilon$, $\exists$ $q$ satisfying
$\dist(p,q)<\epsilon$ such that
\[
F_q^{-1}(1-\delta')> f_p(F_p^{-1}(1-\delta')).
\]}

\ignore{Recall that $\eta$ will be the probability with which some $q\in\cP^{\infty}$
is bankrupted.}

Let $\Phi$ be any loss domination scheme. Recall that $\Phi$
enters on $p$ with probability 1, 
in the sense that the loss dominants set by $\Phi$ will eventually become finite with 
probability $1$ under $p$.
For all
$n \ge 1$, let
\[
\R_{n} :=\sets{x^n: \Phi(x^n)<\infty}
\]
be the set of sequences of length $n$ on which $\Phi$ has entered and
let $N \ge 1$ be a number such that
\begin{equation}
\label{eq:N}
p(\R_N)>1-\alpha/2.
\end{equation}
For any sequence $x^n$, let $A(x^n)$ be the set of symbols that appear
in it. Recall that the head of the distribution $p$, $H_{p,\gamma}$, was
defined in Section~\ref{s:cdf} to be the set $\sets{y\in \cA:y\le 2F_p^{-1}(1-\gamma/2)}$,
where $\cA$ is the support of $p$.
Further, define for all $\gamma>0$
\[
\R_{p,\gamma,n} :=\sets{x^n\in\R_n: A(x^n)\subseteq H_{p,\gamma})}.
\]

\ignore{We will prove that $\exists q\in\cP^{\infty}$ such that 
\begin{equation}
\label{eq:toprove}
q(\text{ $\Phi$ goes bankrupt })\ge (1-\frac2N-2h(\alpha))\Paren{1-\frac1e}.
\end{equation}
By picking $N$ large enough, we can make this probability of error
arbitrarily close to $(1-2h(\alpha))\Paren{1-\frac1e}$.}

Set\footnote{Please note that in the interest of simplicity, 
we have not attempted to provide the best scaling for $\epsilon$ or the tightest
possible bounds in arguments below} $\epsilon=\frac1{16(\ln 2) N^8}$.
Since $p$ is deceptive, there exists $\delta>0$ such that
for all $f(\delta)\in\reals$, there exists a distribution $q\in\cP$ satisfying
both
\begin{equation}
\label{eq:dcpt}
\dist(p,q)<\epsilon=\frac1{16(\ln2) N^8}
\text{ and }
F_{q}^{-1}(1-\delta)> f(\delta).
\end{equation}
While the number
$f(\delta)$ can be arbitrary above, we focus on a specific number
dependent only on $\Phi$.  To define this number, first pick $\gamma_p$
so small that
\begin{equation}
\label{eq:k}
(1-\gamma_p)^{N+1/\delta}\ge
1-\alpha/2.
\end{equation}
Note that the limit of the left side above as $\gamma_p\to0$ is 1, so there is always
some choice $\gamma_p$ that works.
Now, for all $0<\delta'<1$, let
\[
f(\delta') :=
\max_{\substack{x^i\in\R_{p,\gamma_p,i}\\ N\le i \le N+\lceil\frac1{\delta'}\rceil}} 
\Phi(x^i).
\]
\ignore{and let
\[
f_p(0)=\sup_{\substack{0< \delta'\le 1\\F_p^{-1}(1-\delta')=0}} 
\max_{\substack{x^i: N\le i\le N+\lceil\frac1{\delta'}\rceil\\ x^i\in\R_{p,\gamma_p,i}}} \Phi(x^i).
\]}

A word about this parameter $\gamma_p$, since it may not be immediately
apparent why this should be defined.  We will effectively ignore the
$\gamma_p$ tail of the distribution $p$, and focus only on
strings in $\R_{p,\gamma_p,i}$, $N\le i \le N+\frac1\delta$.  The
advantage of doing so is technical---we will be able to handle $p$ and
$q$ as though they were distributions with finite span. This is
crucial since we want to have a finite set over which we take the
supremum on the right side above, so that the maximum is guaranteed to
yield $f(\delta')<\infty$.  Furthermore, note that for $N\le i<
N+\frac1\delta$,
\[
p(\R_{p,\gamma_p,i})\ge 1-\alpha
\] 
from a union bound on~\eqref{eq:N} and~\eqref{eq:k}.

Let $q\in\cP$ satisfy~\eqref{eq:dcpt} with $f(\delta)$ as defined
above.
Applying Lemma~\ref{lm:jn} to distributions over length-$N$
sequences induced by the measures $p,q\in\cP^{\infty}$ corresponding
to the distributions
above,
\[
q(\R_{p,\gamma_p,N})\ge 1-\alpha-\frac2{N},
\]
namely, $\Phi$ has entered with probability (under $q$) at least
$1-\alpha-\frac2{N}$ for length $N$ sequences. Since the insurer
cannot quit once it has entered, the scheme has entered with
probability (under $q$) at least $1-\alpha-\frac2{N}$ for all $n$
length sequences where $n \ge N$. Namely for all $n\ge N$,
\[
q(\R_{p,\gamma_p,n})\ge 1-\alpha-\frac2{N}.
\]

For convenience, let $M=\lceil\frac1\delta\rceil$.
Let the
distribution $q$ be in force.  
We have set things up so that $\Phi$ is bankrupted whenever any element 
in the $\delta$-tail of $q$ follows any sequence
in $R_{p,\gamma_p,i}$, where $N\le i\le N+M-1$. To see this, note that
for all $x\in T_{q,\delta}$
\begin{align}
\label{eq:fq}
x\ge F_q^{-1}(1-\delta)
\ge
f(\delta)
&= 
\max_{\substack{X^i\in\R_{p,\gamma_p,i}\\ N\le i \le N+\lceil\frac1\delta\rceil}} 
\Phi(X^i).
\end{align}
Equivalently, conditioned on any sequence in $\R_{p,\gamma_p,i}$ with
$i$ between $N$ and $N+M-1$, using~\eqref{eq:tail} the scheme $\Phi$
fails with probability (under $q$) at least $\delta$ in step $i+1$.

A sequence on which $\Phi$ has entered, but such that $\Phi$ has not
been bankrupted on any of the sequence's prefixes is called a \emph{surviving}
sequence. 

Consider a surviving sequence $x^N \in\R_{p,\gamma_p,N}$. Given $x^N$, let the conditional
probability that $\Phi$ is bankrupted in the following step be
$\delta_N$. From~\eqref{eq:fq}, as mentioned before, we have $\delta_N\ge
\delta$.
\ignore{ and as mentioned in the note
  before the Theorem, conventions adopted in the definitions of
  cumulative distribution functions in Section~\ref{s:cdf} and
  Equation~\eqref{eq:fq}}

Now, given $x^N \in\R_{p,\gamma_p,N}$, the conditional probability that $\Phi$ is bankrupted
in at most two further steps is,
\[
\delta_N+(1-\delta_N)\delta_{N+1}\ge \delta+(1-\delta)\delta,
\]
where $\delta_{N+1}$ is interpreted as the weighted average (over
surviving length-$(N+1)$ suffixes of $x^N$) of the conditional probability that
$\Phi$ goes bankrupt in step $N+2$ given a surviving sequence of length $N+1$. 

Similarly, given a sequence $x^N \in\R_{p,\gamma_p,N}$, the
probability that $\Phi$ is bankrupted on suffixes of $x^N$
with length between $N$ and $N+M$ is
\[
\delta_N+
(1-\delta_N)\delta_{N+1}
+\ldots+
\delta_{N+M}\prod_{i=N}^{N+M-1}(1-\delta_i)
\]
for some $\delta_N,\delta_{N+1}\upto \delta_{N+M}$, all of which 
are $\ge\delta$.

Let $q_1$ be the probability (under $q$) of all survivors in $\R_{p,\gamma_p,N}$, and
$q_2$ be the probability (under $q$) of all sequences in $\R_{p,\gamma_p,N}$ where $\Phi$
has already been bankrupted.  Therefore $q_1+q_2=q(\R_{p,\gamma_p,N})$.
\ignore{
Therefore $\Phi$ is bankrupted with probability
\begin{align*}
&\ge q_2 +q_1\Paren{\delta_N
+
(1-\delta_N)\delta_{N+1}
+\ldots+
\delta_M\prod_{i=N}^{M}(1-\delta_i)}
\end{align*}
where as before $M=\lceil\frac1{\delta}\rceil$. }

Let $\overline{\delta}$ stand for $1-\delta$. 
Now $\Phi$ is bankrupted with probability 
\begin{align}
\nonumber
&\ge
q_2+
q_1\Paren{\delta_N
+\ldots+
\delta_{N+M}\prod_{i=N}^{N+M-1}(1-\delta_i)}\\
\nonumber
&=
q_2+
q_1
\Paren{\delta_{{}_{N}}
+\overline{\delta_{{}_N}}\Paren{
   \delta_{{}_{N+1}}+
   \overline{\delta_{{}_{N+1}}}\Paren{
\ldots \Paren{ \delta_{{}_{N+M-1}}+\overline{\delta_{{}_{N+M-1}}}\delta_{{}_{N+M}}}}}}\\
\nonumber
&\ge q_2+
q_1
\Paren{
\delta
+(1-\delta)\delta
+\ldots +
(1-\delta)^{M}\delta}\\
\nonumber
&=
q_2+
q_1
\Paren{1-(1-\delta)^{\lceil 1/\delta \rceil}}\\
\nonumber
&\ge
q(\R_{p,\gamma_p,N})
\Paren{1-(1-\delta)^{\lceil 1/\delta \rceil}}\\
\nonumber
&\ge
\Paren{
1-\alpha-\frac2{N}}  
\Paren{1-(1-\delta)^{\lceil 1/\delta \rceil}}.
\end{align}
The Theorem follows.
\eTheorem

\ignore{ We will show that if $q$ were the distribution in force, it
  would be bankrupted with probability $1-\alpha-\frac1e$.  Let
\[
\Ab{\beta}=\sets{ y\in\naturals : y \le F_p^{-1}(1-\beta)}
\]
be the set of all elements in the support of $p$ except a tail of at most
$\beta$, and $\Ab{\beta}^n$ is the set of $n$ element sequences each of whose
symbols is taken from $\Ab{\beta}$. 
Consider the set of \emph{relevant} sequences,
\[
R=\Ab{\frac1{n^{*}}}^{n^{*}} \cap \hat\varphi(n^{*}).
\]
Since $p(R)\ge 1-\alpha-\frac1e$ and $D(p^{n^{*}}||q^{n^{*}})\le \alpha$, it follows
from Lemma~\ref{lm:app} that 
\[
q(R)\ge 1-\alpha-h\Paren{1-\alpha-\frac1e}.
\]
Note that for $\alpha$ small enough, $q(R)>0$. For example, for
$\alpha=.02$, $q(R)=.0166$.
Consider a sequence in $R$.  From~\eqref{eq:fq}, the scheme $\Phi$
gets bankrupted within the next $n^{*}$ steps if any element in the
$\frac1{n*}$ tail of $q$ occurs within the next $n^{*}$
steps. Conditioned on $R$, the probability (under $q$) of this
happening is $1-\Paren{1-\frac1{n^{*}}}^{n^{*}}$. 
Therefore, the probability under $q$ of $\Phi$ being bankrupted is at least,
\[
\Paren{1-\frac1e}q(R). \eqed
\]
\eTheoremp}

\ignore{\bCorollary
The set of \iid monotone distributions is not insurable.
\eCorollary}
\section{Sufficient condition for insurability}
\label{s:sff}
When no $p\in\cP$ is deceptive, given any $\eta>0$ we will construct a loss domination scheme that 
goes bankrupt with probability $\le\eta$.

\ignore{ The necessary condition in Section~\ref{s:ncs} is also sufficient for insurability. With the convention for
  cumulative distribution functions adopted in Section~\ref{s:cdf}, recall that when the marginal loss distribution is
  $q$ the probability that a loss domination scheme is bankrupted by the loss incurred in period $i$ is strictly less
  than $\delta$ if the loss-dominant set at the beginning of that period is at least $2F_q^{-1}(1-\delta/2)$.  Given a
  sequence of losses $x^n$ we will keep track of its empirical distribution (which we call its \emph{type}) in a
  somewhat unconventional, unnormalized way.  The type of $x^n$ will be written as a sequence of unnormalized fractions
  indexed by integers with common denominator $n$ (we only write the nontrivial part of the sequence).  Thus the string
  123111, which is of length $6$, has type $(4/6,1/6,1/6)$, while the string $24224421$, of length $8$, has type
  $(1/8,4/8,0/8,3/8)$.  }

If no $p\in \cP$ is deceptive, there is for each $p \in \cP$ a number
$\epsilon_p>0$ such that, for every percentile $\delta > 0$, 
there is a uniform bound on the $\delta$-percentile over the set of
probability distributions in the neighborhood
\[
\sets{ p'\in \cP: \dist(p',p)< \epsilon_p},.
\]
We pick such an $\epsilon_p$ for each $p \in \cP$ and call it
the \emph{reach} of $p$.
For $p\in\cP$, the set
\[
B_p = \sets{ p'\in\cP: \dist(p,p')  < \epsilon_p},
\]
where $\epsilon_p$ is the reach of $p$, will play the role of the set
of probability distributions in $\cP$ for which it
will be okay to eventually set loss-dominants assuming $p$ is in force.

To prove that $\cP^\infty$ is insurable if no distribution among its one dimensional marginals $\cP$ is deceptive, we
will need to find a way to cover $\cP$ with countably many sets of the form $B_p$ above.  Unfortunately, $\dist(p,q)$ is
not a metric, so it is not immediately clear how to go about doing this. On the other hand note that $\dist(p',p)\le
|p-p'|_1/\ln 2$, where $|p-p'|_1$ denotes the $\ell_1$ distance between $p$ and $p'$ (see Lemma \ref{lm:dist} in the
Appendix). Therefore, we can instead bootstrap off an understanding of the topology induced on $\cP$ by the $\ell_1$
metric.

\subsection{Topology of $\cP$ with the $\ell_1$ metric}\label{ss:lindel}
The topology induced
on $\cP$ by the $\ell_1$ metric is
Lindel\"of, i.e. any covering of $\cP$ with open sets in the $\ell_1$ topology
has a countable subcover (see~\cite[Defn. 6.4]{Dugundji}
for definitions and properties of Lindel\"of topological spaces).

We can show that $\cP$ with the $\ell_1$ topology is Lindel\"of by
appealing to the fact that the set of all probability distributions
on $\naturals$ with the $\ell_1$ topology, is second countable, i.e. 
that it has a countable basis. The set of all distributions
on $\naturals$ along with $\ell_1$ topology has a countable basis 
because it has a countable norm-dense set (consider the set of all 
probability distributions on $\naturals$ with finite support and with 
all probablities being rational). Now, $\cP$, as a topological 
subspace of a second countable topological space is also second 
countable~\cite[Theorem 6.2(2)]{Dugundji}. Finally, every second countable 
topological space is Lindel\"of \cite[Thm. 6.3]{Dugundji}, hence
$\cP$ is Lindel\"of.

\subsection{Sufficient condition}
We now have the machinery required to prove that if no $p\in\cP$ is
deceptive, then $\cP^\infty$ is insurable, which is the other direction of 
Theorem \ref{thm:ncssff}, as stated next.

\bTheorem\label{thm:sff} 
If no $p\in\cP$ is deceptive, then $\cP^{\infty}$ is insurable.
\Proof
The proof is constructive. For any $0<\eta<1$, we obtain a loss domination scheme $\Phi$
such that for all $p\in\cP^\infty$, $p\Paren{\Phi \text{ goes bankrupt }}<\eta$.

For $p\in\cP$, let
\[
Q_p = \Sets{ q: |p-q|_1 < \frac{{\epsilon_p}^2(\ln 2)^2}{16}},
\]
where $\epsilon_p$ is the reach of $p$. We will call $Q_p$ as the \emph{zone} of $p$.
The set $Q_p$ is non-empty when $\epsilon_p>0$.

For large enough $n$, the set of loss sequences of length $n$ with
empirical distribution in $Q_p$ will ensure that the loss domination
scheme $\Phi$ to be proposed enters with probability 1 when $p$ is in
force.  Note that if $\epsilon_p>0$ is small enough then $Q_p\cap \cP
\subset B_p$---we will assume wolog that $\epsilon_p>0$ is always taken
so that $Q_p\cap \cP \subset B_p$.

Since no $p\in\cP$ is deceptive, none of the zones $Q_p$ are empty and
the space $\cP$ of distributions can be covered by the sets $Q_p\cap
\cP$, namely
\[
\cP=\union_{p\in\cP} (Q_p\cap \cP).
\] 
From Section \ref{ss:lindel}, we know that $\cP$ is Lindel\"of under
the $\ell_1$ topology. Thus, there is a countable set
$\tilde\cP\subseteq\cP$, such that $\cP$ is covered by the collection
of relatively open sets 
\[
\sets{Q_{\tilde p}\cap \cP: {\tilde p}\in\tilde\cP}.
\]
We let the above collection be denoted by $\cQ_{\tilde\cP}$.
We will refer to $\tilde\cP$ as the \emph{quantization} of $\cP$ and to elements of $\tilde\cP$ 
as \emph{centroids} of the quantization, borrowing from commonly used literature in classification.

\ignore{In an abuse of notation, let $2^\naturals$ be the set of all
  finite subsets of $\naturals$. Note that $2^\naturals$ is
  countable. To see this, represent every number by its Elias code,
  and every finite subset of naturalsn by concatenating the
  representations of its components. Each subset of $\naturals$ hence
  maps to a unique rational number. Since rational numbers are
  countable, so is $\naturals^*$.}  

We index the countable set of centroids, $\tilde\cP$ (and reuse the
index for the corresponding elements of $\cQ_{\tilde \cP}$) by
$\iota:{\tilde \cP}\to \naturals$.

We now describe the loss domination scheme $\Phi$ having the property that
for all $p\in\cP^\infty$, 
\[
p\Paren{\Phi \text{ goes bankrupt }}<\eta.
\]
\paragraph{Preliminaries} Consider a length-$n$ sequence $x^n$
on which $\Phi$ has not entered thus far. Let the empirical distribution of the sequence
be $q$, and let
\[
\cP_q' := \sets{ p'\in{\tilde \cP}: q\in Q_{p'} }
\]
be the set of centroids in the quantization of $\cP$ (elements of $\tilde\cP$) which can potentially \emph{capture} $q$.
Note that $q$ in general need not belong to $\tilde{\cP}$ or $\cP$.

If $\cP_q'\ne \emptyset$, we will further refine the set of
distributions that could capture $q$ further to $\cP_q \subset \cP_q'$
as described below.  Refining $\cP_q'$ to $\cP_q$ ensures that
models in $\cP_q'$ do not prematurely capture loss sequences.

Let $p$ be the model in force, which remains unknown. The idea is that
we want sequences generated by (unknown) $p$ to be captured by those centroids
of the quantization $\tilde\cP$ that have $p$ in their reach. We will
require~\eqref{eq:bnkrpt} below to ensure that the probability (under
the unknown $p$) of all sequences that may get captured by 
centroids $p'\in\cP_q$ not having $p$ in its reach remains small.  In
addition, we impose~\eqref{eq:bnkrpttwo} as well to resolve a
technical issue since $q$ need not, in general, belong to $\cP$.

For $p'\in\cP_q'$, let the reach of $p'$ be
$\eps$, and define
\[
D_{{}_{p'}}
:=
\frac{{\epsilon_{p'}}^4 (\ln 2)^4}{256}~.
\]
In case the underlying distribution $p$ happens to be out of the reach
of $p'$ (wrong capture), the quantity $D_{{}_{p'}}$ will later lower
bound the distance of the empirical $q$ in question from the
underlying $p$.

Specifically, we place $p'$ in $\cP_q$ if $n$ satisfies
\begin{equation}
\label{eq:bnkrpt}
\exp\Paren{-nD_{{}_{p'}}/18}
\le 
\frac{\eta }{2C(p') \iota(p')^2 n(n+1)},
\end{equation}
and 
\begin{equation}
\label{eq:bnkrpttwo}
2F_q^{-1}(1-\sqrt{D_{{}_{p'}}}/6) \le \log C(p'),
\end{equation}
where $C(p')$ is 
\[
C(p') := 2^{2\Paren{{\sup_{r\in B_{p'}}} F_r^{-1}(1-\sqrt{D_{{}_{p'}}}/6)}}.
\]
Note that $C(p')$ is finite since $p'$ is not deceptive. Comparison
with Lemma~\ref{lm:yeung} will give a hint as to why the equations
above look the way they do.

\paragraph{Description of $\Phi$}
For the sequence $x^n$ with type $q$, if $\cP_q=\emptyset$, the scheme 
does not enter yet. If $\cP_q\ne\emptyset$, let $p_q$ denote the distribution 
in $\cP_q$ with the smallest index.

All sequences with prefix $x^n$ (namely sequences obtained by concatenating 
$x^n$ with by any other sequence of symbols) are then said to be \emph{trapped} 
by $p_q$---namely, loss-dominants will be based on $p_q$. The loss-dominant assigned
for a length-$m$ sequence trapped by $p_q$ is
\[
2g_{p_q}\Paren{\frac{\eta}{4n(n+1)}}
:=
2\sup_{r\in B_{p_{{}_q}}} F^{-1}_r\Paren{1-\frac{\eta}{4n(n+1)}}.
\]

\paragraph{$\Phi$ enters with probability 1} 
First, we verify that the scheme enters with probability 1, no matter what distribution $p\in\cP$ is in force. Every
distribution $p\in\cP$ is contained in at least one of the elements of the cover $\cQ_{\tilde \cP}$. 

Recall the enumeration of $\tilde\cP$. Let $p'$ 
be centroid with the smallest index among all centroids
in ${\tilde \cP}$ whose zones contain $p$.
Let $Q$ be the zone of $p'$.
There is thus some $\gamma > 0$ such that the neighborhood around $p$ given
by
\[
I(p,\gamma) := \sets{q: |p-q|_1 <\gamma}
\]
satisfies $I(p,\gamma)\subseteq Q$.
Note in particular that $p$ is in the reach of $p'$.

With probability 1, sequences generated by $p$ will have their empirical distribution within $I(p,\gamma)$
(see~\cite{Chu61} or Lemma~\ref{lm:yeung} for an alternate proof). Next~\eqref{eq:bnkrpt} will hold for all sequences
whose empirical distributions that fall in $I(p,\gamma)$ whose length $n$ is large enough---since $C(p')$ and
$\iota(p')$ do not change with $n$, the right hand side diminishes to zero polynomially with $n$ while the left hand
side diminishes exponentially to zero. Thus we conclude~\eqref{eq:bnkrpt} will be satisfied with probability 1.

Next,~\eqref{eq:bnkrpttwo} will also hold almost surely, for if $q$ is the empirical probability of sequences generated
by $p$, then (with a little abuse of notation)
\[
F_q^{-1}(1-\sqrt{D_{{}_{p'}}}/6)\to F_p^{-1}(1-\sqrt{D_{{}_{p'}}}/6)
\] 
with probability 1. Note that the quantity on the left is actually a random variable that is sequence dependent (since $q$ 
is the empirical distribution of the sequence). Furthermore, we also have
\begin{align*}
2F_p^{-1}(1-\sqrt{D_{{}_{p'}}}/6)
&\le2\Paren{{\sup_{r\in B_{p'}}} F_r^{-1}(1-\sqrt{D_{{}_{p'}}}/6)}\\
&=
\log C(p'),
\end{align*}
where the first inequality follows since $p$ is in the reach of $p'$.

Thus the scheme enters with probability 1 no matter which $p\in\cP$ is in force.

\paragraph{Probability of bankruptcy $\le \eta$}
We now analyze the scheme. Consider any $p\in\cP$.  Among sequences on
which $\Phi$ has entered, we will distinguish between those that are
in \emph{good} traps and those in \emph{bad} traps. If a sequence
$x^n$ is trapped by $p'$ such that $p\in B_{p'}$, $p'$ is a good
trap. Conversely, if $p\notin B_{p'}$, $p'$ is a bad trap.

\textit{(Good traps)}
Suppose a length-$n$ sequence $x^n$ is in a good trap, namely, it is
trapped by a distribution $p'$ such that $p\in B_{p'}$. Recall that
the loss-dominant assigned is
\[
2g_{p'}\Paren{\frac{\eta}{4n(n+1)}}
\ge 
2F_p^{-1}\Paren{1-\frac{\eta}{4n(n+1)}},
\]
where the inequality follows because $p'$ is not deceptive, and $p$ is
within the reach of $p'$. Therefore from~\eqref{eq:head}, given any
sequence in a good trap the scheme is bankrupted with conditional
probability at most $\delta'=\eta/2n(n+1)$ in the next
step. Therefore, summing over all $n$, sequences in good traps
contribute at most $\eta/2$ to the probability of bankruptcy.

\textit{(Bad traps) } We will show that the probability with which
sequences generated by $p$ fall into bad traps
$\le\eta/2$. Pessimistically, the conditional probability of
bankruptcy in the very next step given a sequence falls into a bad
trap is going to be upper bounded by 1. Thus the contribution to
bankruptcy by sequences in bad traps is at most $\eta/2$.

Let $q$ be any length-$n$ empirical distribution trapped by
$\tilde{p}$ with reach $\tilde\epsilon$ such that $p\notin
B_{\tilde{p}}$. 

If $p$ is ``far'' from $\tilde{p}$ (because $p$ is not in
$\tilde{p}$'s reach), namely
\[
\dist(\tilde{p},p)\ge \tilde\epsilon,
\]
but $q$ is ``close'' to $\tilde{p}$ (because $q$
has to be in $\tilde{p}$'s zone to be captured by it), namely
\[
|\tilde{p}-q|_1 
<
\frac{{\tilde\epsilon}^2(\ln 2)^2}{16},
\]
then we would like $q$ to be far from $p$. That is exactly what we obtain
from the triangle-inequality like Lemma~\ref{lm:dpq}, namely that
\[
\dist(p,q)\ge\frac{\tilde\epsilon^2\ln 2}{16}
\] 
and hence, for all $q$
trapped by $\tilde{p}$ that
\[
|p-q|_1^2 
\ge \dist^2(p,q)(\ln 2)^2
\ge\frac{\tilde{\epsilon}^4 (\ln 2)^4}{256}
=D_{{}_{\tilde{p}}}^2.
\]
We need not be concerned that the right side above depends on $\tilde{p}$,
and there may be actually no way to lower bound the rhs as a function of just $p$. Rather,
we take care of this issue by setting the entry point appropriately via~\eqref{eq:bnkrpt}.

Thus, for $p\in\cP^\infty$, the probability length-$n$ sequences with empirical distribution 
$q$ is trapped by a bad $\tilde{p}$ is, using~\eqref{eq:bnkrpt} and~\eqref{eq:bnkrpttwo}
\begin{align*}
&\le
p\biggl( 
| q-p|^2 \ge D_{{}_{\tilde{p}}}
\text{ and } 
\ignore{\biggr.\\
&\qquad\qquad\qquad }
2F_q^{-1} (
1-\frac{\sqrt{D_{{}_{\tilde{p}}}}}6
)
\le \log C(\tilde{p}) 
\biggr)\\
\ignore{&=
p\biggl(
| q-p| \ge \sqrt{D_{{}_{\tilde{p}}}}
\text{ and } \biggr.\\
&\qquad\qquad\qquad F_q^{-1}(1-\sqrt{D_{{}_{\tilde{p}}}}/3) \le \log C(\tilde{p})
\biggr)\\}
&\ale{(a)}
(C(\tilde{p})-2) \exp\Paren{-\frac{n D_{{}_{\tilde{p}}}}{18}}\\
&\ale{(b)}
\frac{\eta (C(\tilde{p})-2)}
{2 C(\tilde{p}) \iota(\tilde{p})^2 n(n+1)} \\
&\le \frac{\eta}{2\iota(\tilde{p})^2 n(n+1)},
\end{align*}
where the inequality $(a)$ follows from Lemma~\ref{lm:yeung}
and $(b)$ from~\eqref{eq:bnkrpt}. 
Therefore, the probability of sequences
falling into bad traps
\[
\le
\sum_{n\ge1}\sum_{\tilde{p}\in {\tilde \cP}}
\frac{\eta}{2\iota(\tilde{p})^2 n(n+1)} 
\le\eta/2
\]
since
$\sum_{\tilde{p}\in {\tilde \cP}}
\frac1{\iota(\tilde{p})^2}
\le
\sum_{n\ge1}
\frac1{n(n+1)}
=
1.
$
The theorem follows.
\eTheorem

\section{Concluding remarks}	\label{s.concrem}

The loss domination problem formulated and solved in this paper appears to be of 
natural interest. However, there are several features of the insurance problem 
formulated here that might appear troubling even to the casual reader. In practice
an insured party entering into an insurance contract would expect some stability in the 
premiums that are expected to be paid. A natural direction for further research is
therefore to study how the notion of insurability of a model class changes when 
one imposes restrictions on how much the premium set by the insurer can vary from
period to period. Another obvious shortcoming of the formulation of the insurance
problem studied here is the assumption that the insured will accept any contract
issued by the insurer. Since the insured in our model represents an aggregate of 
individual insured parties, a natural direction to make the framework more realistic
would be to think of the insured parties as being of different \emph{types}. This would
in effect make the total realized premium from the insured (the aggregate of the insured
parties) and the distribution of the realized loss in each period a function of the size of the premium per insured party set by the insurer in that period. 
Characterizing which model classes are insurable when the 
realized premium and the realized loss are functions of a set premium per insured party 
would be of considerable interest.

Both for the loss domination problem and for the insurance problem, working with model
classes for the loss sequence that allow for dependencies in the loss from period to period,
for instance Markovian dependencies, would be another interesting direction for further
research. Considering models with multiple, possibly competing insurers, as well as considering
an insurer operating in multiple markets, where losses in one market can be offset
by gains in another, also seem to be useful directions to investigate.

\section*{Acknowledgments}
We thank C. Nair (Chinese Univ of Hong Kong) and K. Viswanathan (HP
Labs) for helpful discussions.  N. Santhanam was supported by NSF
Grants CCF-1065632, CCF-1018984 and EECS-1029081.  V. Anantharam was
supported by the ARO MURI grant W911NF- 08-1-0233, “Tools for the
Analysis and Design of Complex Multi-Scale Networks”, the NSF grant
CNS-0910702, the NSF Science \& Technology Center grant CCF-0939370,
“Science of Information”, Marvell Semiconductor Inc., and the
U.C. Discovery program.

\section*{Appendix}
\ignore{
\bLemma
\label{lm:app} 
Let $p$ and $q$ be probability distributions on a countable
set $\cA$. Suppose
$\dist(p,q)\le\epsilon$.  For any $S\subset \cA$ and $\alpha<1-\ln
2=.30685$, if $p(S)\ge 1-\alpha$, then
\[
q(S)\ge 1-2\epsilon-2h(\alpha).
\]
\Proof 
Let $p_S$ (respectively $q_S$) denote a binary
distribution, the two probabilities of which correspond to
$[p(S),1-p(S)]$ (respectively $[q(S),1-q(S)]$). Now,
\begin{align*}
\epsilon
\ge \dist(p,q)
&\ge D\Paren{p_S||\frac{q_S+p_S}2}\\
&\ge
p(S)\log\frac2{p(S)+q(S)}
-
h(p_S)\\
&\ge
(1-\alpha)\log\frac2{p(S)+q(S)}
-
h(\alpha).
\end{align*}
The last inequality follows because the condition $p(S)\ge1-\alpha\ge\half$ 
implies $h(p_S)\le h(1-\alpha)$.
Therefore,
\[
\log\frac2{p(S)+q(S)} \le \frac{h(\alpha)+\epsilon}{(1-\alpha)},
\]
implying that 
\begin{align*}
\frac{1+q(S)}2 \ge \frac{p(S)+q(S)}2&\ge 2^{-\frac{h(\alpha)+\epsilon}{1-\alpha}}\\
&\ge 1-\Paren{h(\alpha)+\epsilon},
\end{align*}
where the last inequality follows since $\ln 2 \le 1-\alpha$. 
\eLemma}

\bLemma
\label{lm:dist}
Let $p$ and $q$ be probability distributions on $\naturals$. Then
\[
\frac1{4\ln 2} |p-q|_1^2 \le \dist(p,q) \le \frac1{\ln2} |p-q|_1~.
\]
If, in addition, $r$ is a probability distribution on $\naturals$,
then 
\[
\dist(p,q)+ \dist(q,r)\ge \dist^2(p,r)\frac{\ln 2}{8}.
\]
\Proof
The lower bound in the first statement follows since
\[
D\Paren{p||\frac{p+q}2}\ge \frac1{2\ln2} \frac14 |p-q|_1^2
\]
and similarly for $D\Paren{q||\frac{p+q}2}$.
Since $\ln(1+z)\le z$ for all $z\ge0$, the upper bound in the first statement follows 
as below:
\begin{align*}
\dist(p,q) \ln 2&\le
\sum_{x: p(x)\ge q(x)} p(x) \Paren{\frac{p(x)-q(x)}{p(x)+q(x)}}
+
\sum_{x': q(x')\ge p(x')} q(x') \Paren{\frac{q(x')-p(x')}{p(x')+q(x')}}\\
&\le
|p-q|_1.
\end{align*}
To prove the triangle-like inequality, note that
\begin{align*}
\dist(p,q)+ \dist(q,r)
&\ge 
\frac1{4\ln2} \Paren{|p-q|^2_1+|q-r|^2_1}\\
&\ge
\frac1{8\ln2} \Paren{|p-q|_1+|q-r|_1}^2\\
&\ge
\frac1{8\ln2} \Paren{|p-r|_1}^2\\
&\ge
\frac{\ln2}{8} \dist(p,r)^2,
\end{align*}
where the last inequality follows from the upper bound on $\dist(p,r)$ 
already proved.
\eLemma

\bLemma
\label{lm:jn}
Let $p$ and $q$ be probability distributions on a countable set $\cA$
with $\dist(p,q)\le \epsilon$.
Let $p^N$ and $q^N$ be distributions over $\cA^N$ obtained by \iid
sampling from $p$ and $q$ respectively (the distribution induced by
the product measure). For any $R_N\subset \cA^N$ and $\alpha>0$, if $p^N(R_N)\ge 1-\alpha$, then
\[
q^N(R_N)\ge 1- \alpha-2N^3\sqrt{4\epsilon\ln 2}-\frac1N.
\]
\Proof
Let 
\[
\cB_1 = \Sets{ i\in\cA : q(i) \le p(i)\Paren{1-\frac1{N^2}} },
\]
and let 
\[
\cB_2 = \Sets{ i\in\cA : p(i) \le q(i)\Paren{1-\frac1{N^2}} },
\]
If $\dist(p,q)\le \epsilon$, then we have
\[
\sqrt{\epsilon} \ge \sqrt{\dist(p,q)} \ge \frac{|p-q|_1}{\sqrt{4\ln2}}.
\]
It can then be easily seen that 
\begin{equation}
\label{eq:tmp}
p(\cB_1\cup \cB_2) \le 2N^2\sqrt{4\epsilon\ln2}
\text{ and }
q(\cB_1\cup \cB_2) \le 2N^2\sqrt{4\epsilon\ln2}
\end{equation}
because
\[
|p-q|_1
\ge
\sum_{x\in\cB_1} (p(x) - q(x) )
\ge
\frac{p(\cB_1)}{N^2}
\ge
\frac{q(\cB_1)}{N^2}
\]
and similarly
\[
N^2 |p-q|_1
\ge
q(\cB_2)
\ge
p(\cB_2).
\]
Let $S=\cA-\cB_1\cup\cB_2$. We have for all $x\in S$,
\begin{equation}
\label{eq:tmptwo}
q(x) \ge p(x) \Paren{1-\frac1{N^2}}.
\end{equation}
and from~\eqref{eq:tmp} we have $p(S)\ge 1-{2N^2}\sqrt{4\epsilon\ln2}$.
\ignore{From Lemma~\ref{lm:app},
\[
q(S) \ge 1-2\epsilon-2N^2\sqrt{4\epsilon\ln2}.
\]}
Now, we focus on the set $S_N\subset \cA^N$ containing all length-$N$
strings of symbols from $S$. Clearly 
\[
p(S_N) \ge 1- 2N^3\sqrt{4\epsilon\ln2}.
\]
Thus we have
\[
p(R_N\cap S_N) \ge 1-2N^3\sqrt{4\epsilon\ln2}-\alpha.
\]
From~\eqref{eq:tmptwo}, for all $x^N\in S_N$,
\[
q(x^N) \ge p(x^N)\Paren{1-\frac1{N^2}}^N\ge p(x^N)\Paren{1-\frac1N}.
\]
Therefore, 
\[
q(R_N) \ge q(R_N\cap S_N) \ge (1-2N^3\sqrt{4\epsilon\ln2}-\alpha)\Paren{1-\frac1N}
\ge 1-\alpha-2N^3\sqrt{4\epsilon\ln2}-\frac1N.\eqed
\]
\eLemmap

\bLemma
\label{lm:dpq} 
Let $\epsilon_0 > 0$. If
\[
|p_0-q|_1\le \frac{\epsilon_0^2(\ln 2)^2}{16}~,
\]
then
for all $p\in\cP$ with $\dist(p,p_0)\ge\epsilon_0$, we have
\[
\dist(p,q)\ge\frac{\epsilon_0^2\ln 2}{16}.
\]
\Proof
Since 
\[
|p_0-q|_1\le \frac{\epsilon_0^2(\ln 2)^2}{16},
\]
Lemma~\ref{lm:dist} implies that
\[
\dist(p_0,q)\le\frac{\epsilon_0^2\ln 2}{16}.
\]
Further,
Lemma~\ref{lm:dist} then implies that
\[
\dist(p,q)+\frac{\epsilon_0^2\ln 2}{16}
\,\ge\,
\dist(p,q)+\dist(p_0,q)
\,\ge\,
\frac{\dist^2(p,p_0)\ln 2}{8}
\,\ge \,
\frac{\epsilon_0^2\ln 2}{8},
\]
where the last inequality follows since $\dist(p,p_0)\ge\epsilon_0$.
\eLemma

\bLemma
\label{lm:yeung}
Let $p$ be any probability distribution on $\naturals$. Let $\delta>0$
and let $k \ge 2$ be an integer.
Let $X_1^n$ be a sequence generated \iid  with marginals $p$ and let $q(X^n)$ be
the empirical distribution of $X_1^n$.
Then 
\begin{align*}
p\Paren{ |q(X^n)-p|>\delta\text{ and }2F_q^{-1}(1-\delta/6) \le k} 
&\le 
(2^{k}-2) \exp\Paren{-\frac{n\delta^2}{18}}.
\end{align*}
\bRemark
There is a lemma that looks somewhat similar in~\cite{HY10}. The difference from~\cite{HY10} is that the right
side of the inequality above does \emph{not} depend on $p$, and this property is crucial for its use here.
\eRemark
\Proof 
The starting point is the following result. Suppose $p'$ is a probability distribution on 
$\naturals$ with finite support of size $L$. 
Then from~\cite{WOSVW05}, if we consider length $n$ sequences,
\begin{equation}\label{eq:base}
p'( |q(X^n)-p'|_1\le t) \ge 1- (2^L-2) \exp\Paren{-\frac{n t^2}{2}}.
\end{equation}
Since $k\ge 2$, consider the distributions 
$p'$ and $q'$ with support $A=\sets{1\upto k-1}\union \sets{-1}$, obtained as
\[
p'(i)=
\begin{cases}
p(i) & 1\le i < k\\
\sum_{j=k}^\infty p(j) & i = -1,
\end{cases}
\]
and similary for $q'$. 

From~\eqref{eq:base},
\[
p'(|p'-q'|_1 > \delta/3)
\le
(2^{k}-2) \exp\Paren{-\frac{n\delta^2}{18}}.
\]
We will see that all sequences generated by $p$ with empirical
distributions $q$ satisfying
\[
|p-q|_1 > \delta
\text{ and }
2F_q^{-1}(1-\delta/6) \le k
\]
are now mapped into sequences generated by $p'$ with empirical $q'$ satisfying
\begin{equation}
\label{eq:toshow}
|p'-q'|_1 > \delta/3
\text{ and }
q'(-1)\le \delta/3.
\end{equation}
Thus, we will have
\begin{align*}
&p( |q(X^n)-p|_1 > \delta \text{ and } 2F_q^{-1}(1-\delta/6) \le k) \\
&\le
p'(|p'-q'|_1 > \delta/3 \text{ and } q'(-1)\le \delta/3) \\
&\le
(2^{k}-2) \exp\Paren{-\frac{n\delta^2}{18}}.
\end{align*}
\ignore{where the first inequality follows by simple logical implications and
by observing that $p(\naturals- \sets{1\upto k(q)-1})=p'(-1)$. The lemma
then follows.}

Finally we observe~\eqref{eq:toshow} as in~\cite{HY10}
\begin{align*}
|p-q|_1 &-\sum_{l=1}^{k-1} |p(l)-q(l)|\\ 
&\le\sum_{j=k}^{\infty} (p(j)-q(j)) +2\sum_{j=k}^{\infty} q(j)\\
&\le|p'(-1)-q'(-1)|+2\delta/3,
\end{align*}
where the last inequality above follows from~\eqref{eq:head}.
Since $p(l)=p'(l)$ and $q(l)=q'(l)$ for all $l=1\upto k-1$,
we have
\[
|p'-q'|_1 \ge |p-q|_1 -2\delta/3.
\]
If $|p-q|_1 \ge \delta$ in addition, $|p'-q'|_1 \ge \delta/3$.
\eLemma
\ignore{
The following lemma proved in the Appendix show that $\dist(p,q)$ is
related to the $\ell_1$ norm and satisfies a modification of the
triangle inequality.
\bLemma
\label{lm:dist} For any two distributions $p,q$,
\[
\frac1{4\ln 2} |p-q|_1^2 \le \dist(p,q) \le \frac1{\ln2} |p-q|_1
\]
and for any three distributions $p,q,r$,
\[
\dist(p,q)+ \dist(q,r)\ge \dist^2(p,r)\frac{\ln 2}{8}.\eqed
\]
\eLemmap}

\bibliography{univcod}
\end{document}